\documentclass[a4paper]{article}
\usepackage{amsthm,amsfonts}
\usepackage{amsmath,amssymb}
\usepackage[margin=1.2in]{geometry}
\usepackage{float}
\usepackage{latexsym}
\usepackage{graphicx}
\usepackage[english]{babel}
\usepackage{mathrsfs}
\usepackage{pgf,tikz}
\usepackage[normalem]{ulem}
\usepackage{tabularx}
\usepackage{enumitem}
\usepackage{ae}
\usepackage[T1]{fontenc}

\newtheorem{Theorem}{Theorem}
\newtheorem{Proposition}{Proposition}

\newtheorem{Remark}{Remark}
\newtheorem{Lemma}{Lemma}
\newtheorem{Example}{Example}
\newtheorem{Definition}{Definition}

\begin{document}

\begin{center}
{\Large\bf Strong difference families of special types\footnote{Partially supported by Italian Ministry of Education under Grant PRIN 2015 D72F16000790001 (S. Costa). Supported by NSFC under Grant 11971053 (Y. Chang), NSFC under Grant 11871095 (T. Feng), NSFC under Grant 11771227, and Zhejiang Provincial Natural Science Foundation of China under Grant LY17A010008 (X. Wang).}}

\vskip12pt

Yanxun Chang$^a$, Simone Costa$^b$, Tao Feng$^a$, Xiaomiao Wang$^c$\\[2ex]
{\footnotesize $^a$Department of Mathematics, Beijing Jiaotong University, Beijing 100044, P. R. China}\\
{\footnotesize $^b$Dipartimento DII, Universit\`a degli Studi di Brescia, Via Branze 38 - 25123 Brescia, Italy}\\
{\footnotesize $^c$Department of Mathematics, Ningbo University, Ningbo 315211, P. R. China}\\
{\footnotesize yxchang@bjtu.edu.cn, simone.costa@unibs.it, tfeng@bjtu.edu.cn, wangxiaomiao@nbu.edu.cn}
\vskip12pt

\end{center}

\noindent {\bf Abstract:}
Strong difference families of special types are introduced to produce new relative difference families from the point of view of both asymptotic existences and concrete examples. As applications, group divisible designs of type $30^u$ with block size $6$ are discussed, $r$-rotational balanced incomplete block designs with block size $6$ are derived for $r\in\{6,10\}$, and several classes of optimal optical orthogonal codes with weight $5$, $6$, $7$, or $8$ are obtained.



\section{Introduction}

Throughout this paper, sets and multisets will be denoted by curly braces $\{\ \}$ and
square brackets $[\ ]$, respectively. Every union will be understood as multiset
union with multiplicities of elements preserved. $A\cup A\cup \dots \cup A$ ($h$ times) will be denoted
by $\underline{h}A$. If $A$ and $B$ are multisets defined on a multiplicative group, then $A\cdot B$ denotes the multiset $[ab:a\in A,b\in B]$. Denote by $\mathbb{Z}_v$ the cyclic group of order $v$.

\begin{Definition}\label{SDF}
Let $\Sigma=[F_1,F_2,\dots,F_s]$ be a family of multisets of size $k$ of a group $(G,+)$ of order $g$, where $F_i=[f_{i,1},f_{i,2},\ldots,f_{i,k}]$ for $1\leq i\leq s$. $\Sigma$ is said to be a $(G,k,\mu)$ {\em strong difference family}, or a $(g,k,\mu)$-SDF over $G$, if the list
$$\Delta \Sigma:=\bigcup_{i=1}^s \Delta F_i:=\bigcup_{i=1}^s [f_{i,a}-f_{i,b}: 1\leq a,b\leq k, a\not=b]=\underline{\mu} G,$$
i.e., every element of $G$ $($0 included$)$ appears exactly $\mu$ times in the multiset $\Delta \Sigma$. The members of $\Sigma$ are called {\em base blocks}.
\end{Definition}

\begin{Example}\label{eg:2,5,20}
$\Sigma=[[0,0,1,1,1],[0,1,1,1,1]]$ is a $(\mathbb{Z}_{2},5,20)$-SDF.
\end{Example}

The concept of strong difference families was introduced in \cite{b99} to provide constructions for relative difference families (cf. also \cite{bf,bz,m}). A $(G,N,k,\lambda)$ {\em relative difference family} (DF), or $(g,n,k,\lambda)$-DF over a group $(G,+)$ of order $g$ relative to a subgroup $N$ of order $n$, is a family $\mathfrak{B}=[B_1,B_2,\dots,B_r]$ of $k$-subsets of $G$, called {\em base blocks}, such that the list
$$\Delta \mathfrak{B}:=\bigcup_{i=1}^r\Delta B_i:=\bigcup_{i=1}^r[x-y:x,y\in B_i, x\not=y]=\underline{\lambda}(G\setminus N),$$
i.e., every element of $G\setminus N$ appears exactly $\lambda$ times in the multiset $\Delta \mathfrak{B}$ and no element of $N$ appears in $\Delta \mathfrak{B}$. When $G$ is cyclic, we say that the $(g,n,k,\lambda)$-DF is {\em cyclic} (cf. \cite{buratti98}).

Write $\mathbb{F}_q$ for the finite field of order $q$ and $\mathbb{F}^*_q$ for its multiplicative group. If $q\equiv 1 \pmod{d}$, then $C_0^{d,q}$ denotes the group of nonzero $d$th powers of $\mathbb{F}_q$ and once a primitive element $\omega$ of $\mathbb{F}_q$ has been fixed, we set $C_i^{d,q}=\omega^i\cdot C_0^{d,q}$ for $i=0,1,\ldots,d-1$.

Inspired by several constructions of relative difference families presented in the literature (see  \cite{abg,b99,b4,b19,bg,byw,ChangYin04,cfw,cfw2,fmi,gy,MaChang04,MaChang05,m,yyl}) which take advantage of suitable ``patterns" of strong difference families, we introduce here the following definition.

\begin{Definition}\label{def}
Let $\Sigma$ be a $(G,k,\mu)$ strong difference family and $d$ be a divisor of $\mu$. $\Sigma$ is said to be {\em of type} $d$ with respect to a prime power $q\equiv 1\pmod{d}$ if there is a collection $\mathfrak{C}$ of subsets of $G\times \mathbb{F}_q$ satisfying the following conditions:
\begin{itemize}
\item[$1)$] the projection of $\mathfrak{C}$ on $G$ coincides with $\Sigma$;
\item[$2)$] $\Delta\mathfrak{C}=\bigcup\limits_{g\in G}\left(\{g\}\times D_g\right)$, where $D_g= C_0^{(q-1)/d,q}\cdot L_g$ and $L_g$ is a multiset of size $\mu/d$ on $\mathbb{F}_q^*$.
\end{itemize}
We will simply say that $\Sigma$ is of type $d$ if this property holds for any prime power $q\equiv 1\pmod{d}$.
\end{Definition}

Let $d'$ be a divisor of $d$ and $q\equiv 1\pmod{d}$ be a prime power, which implies $q\equiv 1\pmod{d'}$. If a strong difference family $\Sigma$ is of type $d$ with respect to $q$, then it is also of type $d'$. This is from the observation that $D_g=C_0^{(q-1)/d,q}\cdot L_g=C_0^{(q-1)/d',q}\cdot S \cdot L_g$, where $S$ is a system of representatives for $C_0^{(q-1)/d',q}$ in $C_0^{(q-1)/d,q}$.

Let $(G,+)$ be a group. For $(x,y)\in G\times \mathbb{F}_q$ and $s\in \mathbb{F}^*_q$, define the notation $(x,y)\cdot(1,s):=(x,ys)$. The following proposition gives an illustration of how Definition \ref{def} is applied to construct relative difference families.

\begin{Proposition} \label{prop:basic}
Let $d$ be a divisor of $\mu$, and $\Sigma$ be a $(G,k,\mu)$ strong difference family of type $d$ with respect to the prime power $q\equiv 1\pmod{\mu}$, which is associated with a collection $\mathfrak{C}$ of subsets of $G\times \mathbb{F}_q$ satisfying the two conditions in Definition $\ref{def}$. If for each $g\in G$, $L_g$ is a system of representatives for $C_0^{\mu/d,q}$ in $\mathbb{F}_q^*$, then given a system of representatives $S$ for $C_0^{(q-1)/d,q}$ in $C_0^{\mu/d,q}$, the family $${\mathfrak B}=[C\cdot \{(1,s)\} : s\in S,\ C\in \mathfrak{C}]$$
is a $(G\times \mathbb{F}_q,G\times \{0\},k,1)$-DF.
\end{Proposition}

\proof One can check that
\begin{align*}
\Delta \mathfrak{B}=&\Delta \mathfrak{C}\cdot \{(1,s):s\in S\}=(\bigcup_{g\in G}\left(\{g\}\times D_g\right))\cdot\{(1,s):s\in S\}\\
=&\bigcup_{g\in G}\left(\{g\}\times (D_g\cdot S)\right)=\bigcup_{g\in G}\left(\{g\}\times (C_0^{(q-1)/d,q}\cdot L_g\cdot S)\right)\\
=&\bigcup_{g\in G}\left(\{g\}\times \mathbb{F}^*_q\right)=G\times \mathbb{F}^*_q.
\end{align*} \qed

\begin{Example}\label{eg:using 2,5,20}
The $(\mathbb{Z}_{2},5,20)$-SDF from Example $\ref{eg:2,5,20}$ is of type $2$ with respect to $q=41$. To prove it, we take $\mathfrak C=\{C_1,C_2\}$, where $C_1=\{(0,0),(0,1),(1,2),(1,8),(1,13)\}$ and $C_2=\{(0,0),(1,3),(1,5)$, $(1,17),(1,25)\}.$ It is readily checked that $D_g=\{1,-1\}\cdot L_g$ for  $g\in \mathbb{Z}_{2}$, where $L_g$ is a system of representatives for $C_0^{10,41}=\{1,-1,9,-9\}$ in $\mathbb{F}_{41}^*$.

Furthermore, take $S=\{1,9\}$ as a system of representatives for $C_0^{20,41}=\{1,-1\}$ in $C_0^{10,41}$. Then by Proposition $\ref{prop:basic}$,
$$\left\{\{(0,0),(0,s),(1,2s),(1,8s),(1,13s)\}, \{(0,0),(1,3s),(1,5s),(1,17s),(1,25s)\} : s\in \{1,9\}\right\}$$
is a $(\mathbb{Z}_{2}\times \mathbb{F}_{41},\mathbb{Z}_{2}\times \{0\},5,1)$-DF.
\end{Example}

To use Proposition \ref{prop:basic}, $C_0^{(q-1)/d,q}$ must be a subgroup of $C_0^{\mu/d,q}$, which implies that $q-1\equiv 0\pmod{\mu}$. That is why we require $q\equiv 1\pmod{\mu}$ in Proposition \ref{prop:basic}.

This paper is devoted to constructing strong difference families of types 2 and 4 (see Propositions \ref{pat2} and \ref{pat4}) in Sections 2 and 3, respectively, by using cyclotomic conditions that generalize the procedure of \cite{cfw,cfw2}. Then apply Proposition \ref{prop:basic} to obtain new asymptotic existence results (see Theorems \ref{ExistenceDF2} and \ref{ExistenceDF4}) and existence results for DFs (see Theorems \ref{DF} and \ref{DF1}). Note that an SDF of type 4 is also of type 2 and, in general, not vice versa. In some circumstances, SDFs of type 4 can produce better asymptotic existence results on DFs than SDFs of type 2 (see Remark \ref{rmk:compare}).

As applications, in Section 4, group divisible designs of type $30^u$ with block size $6$ are discussed, $r$-rotational balanced incomplete block designs with block size $6$ are derived for $r\in\{6,10\}$, and several classes of optimal optical orthogonal codes with weight $5$, $6$, $7$, or $8$ are obtained. We conclude in Section 5 and provide related remarks.

In what follows, we always assume that
\begin{eqnarray*}\label{Q(e,m)}
Q(d,m)=\frac{1}{4}\left(U+\sqrt{U^2+4d^{m-1}m}\right)^2, \mbox{ where } U=\sum_{h=1}^m {m \choose h}(d-1)^h(h-1)
\end{eqnarray*}
for given positive integers $d$ and $m$.
Note that $Q(e,m)<Q(e,m')$ for $m<m'$ and $Q(e,m)<Q(e',m)$ for $e<e'$.
The following theorem ensures the existences of elements satisfying certain cyclotomic conditions in a finite field, and will play a crucial role in proving Theorems \ref{ExistenceDF2} and \ref{ExistenceDF4}.

\begin{Theorem}\label{thm:cyclot bound} {\rm \cite{bp,cj}}
Let $q\equiv 1 \pmod{d}$ be a prime power. Let $\{b_0,b_1,\dots, b_{m-1}\}$ be an arbitrary m-subset of $\mathbb{F}_q$ and $(\beta_0,\beta_1,\dots,\beta_{m-1})$ be an arbitrary element of $\mathbb{Z}_d^m$. Set $X=\{x\in \mathbb{F}_q: x-b_i\in C^{d,q}_{\beta_i} \mbox{ for } i=0,1,\dots,m-1\}$. Then $X$ is not empty for any prime power $q\equiv 1 \pmod{d}$ and $q>Q(d,m)$.
\end{Theorem}

\section{Strong difference families of type 2}

In this section, we shall present a class of strong difference families of type 2, and then use them to establish new (asymptotic) existence results for DFs.

\begin{Proposition}\label{pat2}
Let $(G,+)$ be an abelian group. Let $\Sigma$ be a $(G,k,\mu)$-SDF which is the union of three families $\Sigma_1,\Sigma_2$ and $\Sigma_3$, each of which could be empty, where
\begin{itemize}
\item[$1)$] for any $A\in\Sigma_1$, $A$ is of the form $$[x_{1},\delta+x_{1},x_{2},\delta+x_{2},\dots,x_{ k/2},\delta+x_{ k/2}]\mbox{ if }k\mbox{ is even, and}$$
$$[0,x_{1},\delta+x_{1},x_{2},\delta+x_{2},\dots,x_{\lfloor k/2\rfloor},\delta+x_{\lfloor k/2\rfloor}]\mbox{ if }k\mbox{ is odd},$$
where $\delta$ is either an involution of $G$ or zero, and $\Delta[x_1,x_2,\dots,x_{\lfloor k/2\rfloor}]$ does not contain involutions and zeros;
\item[$2)$] for any $A\in \Sigma_2$, each element of $\Delta(A)$ is either an involution or zero;
\item[$3)$] for any $A\in \Sigma_3$, the multiplicity of $A$ in $\Sigma_3$ is even and $\Delta(A)$ does not contain involutions and zeros.
\end{itemize}
Then $\Sigma$ is of type $2$.
\end{Proposition}
\proof
Let $\Sigma=\Sigma_1\cup\Sigma_2\cup\Sigma_3=[A_1,A_2,\dots,A_n]$, where $\Sigma_3=[A_{t+1},A_{t+2},\dots,A_{n}]$ with $A_{t+1}=A_{t+2},\dots$ and $A_{n-1}=A_n$.
Now, consider the family $\mathfrak{C}=[C_1,C_2,\dots, C_n]$, a collection of subsets of $G\times \mathbb{F}_q$, satisfying that
\begin{itemize}
\item if $k$ is even and $A_i=[x_{i,1},\delta+x_{i,1},x_{i,2},\delta+x_{i,2},\dots, x_{i,k/2},\delta+x_{i,k/2}]\in \Sigma_1$, then $$C_i=\{(x_{i,1},y_{i,1}),(\delta+x_{i,1},-y_{i,1}),\dots, (x_{i,k/2},y_{i,k/2}),(\delta+x_{i,k/2},-y_{i,k/2})\};$$
    if $k$ is odd and $A_i=[0,x_{1},\delta+x_{1},x_{2},\delta+x_{2},\dots,x_{\lfloor k/2\rfloor},\delta+x_{\lfloor k/2\rfloor}]\in \Sigma_1$, then
    $$C_i=\{(0,0),(x_{i,1},y_{i,1}),(\delta+x_{i,1},-y_{i,1}),\dots, (x_{i,\lfloor k/2\rfloor},y_{i,\lfloor k/2\rfloor}),(\delta+x_{i,\lfloor k/2\rfloor},-y_{i,\lfloor k/2\rfloor})\};$$
\item if $A_i=[x_{i,1},x_{i,2},\dots, x_{i,k}] \in \Sigma_2$, then
 $$C_i=\{(x_{i,1},y_{i,1}),(x_{i,2},y_{i,2}),\dots, (x_{i,k},y_{i,k})\};$$
\item if $A_i=A_{i+1}=[x_{i,1},x_{i,2},\dots, x_{i,k}] \in \Sigma_3$ for $i\equiv t+1 \pmod{2}$, then
 $$C_i=\{(x_{i,1}, y_{i,1}),(x_2,y_{i,2}),\dots, (x_{i,k},y_{i,k})\}$$
 and $C_{i+1}= C_i\cdot\{(1,-1)\}.$
\end{itemize}
Then one can check that
$$\Delta {\mathfrak C}=\bigcup_{g\in G} (\{g\}\times D_g),$$
where $D_g=\{1,-1\}\cdot L_g$ and $L_g$ is a multiset of size $\mu/2$ on $\mathbb{F}_q^*$. Note that $\mu$ is necessarily even since the zero element of $G$ clearly appears an even number of times in the list of differences of any multisubset of $G$. It follows that $\Sigma$ is of type 2 with respect to any odd prime power $q$.
\endproof

\begin{Example}\label{eg:2,5,20-con}
Example $\ref{eg:using 2,5,20}$ shows that the $(\mathbb{Z}_{2},5,20)$-SDF $\Sigma=[[0,0,1,1,1],[0,1,1,1,1]]$ from Example $\ref{eg:2,5,20}$ is of type $2$ with respect to $q=41$. Actually $\Sigma=\Sigma_2$ satisfies the hypothesis in Proposition $\ref{pat2}$. Therefore, this SDF is of type $2$ with respect to any odd prime power $q$.
\end{Example}

\begin{Example}\label{10-5-12}
By Proposition $\ref{pat2}$, the following $(\mathbb{Z}_{10},5,12)$-SDF is of type $2$:
\begin{center}
\begin{tabular}{|c|c|c|}
\hline $\Sigma_1$ &$\Sigma_2$ &$\Sigma_3$ \\\hline
$[[0,3,3,7,7]]$&
$[[0, 0, 0, 5, 5]]$&
$[[0, 9, 6, 7, 8],\
[0, 9, 6, 7, 8],\
[0, 8, 4, 6, 7],\
[0, 8, 4, 6, 7]]$\\
\hline
\end{tabular} .
\end{center}
\end{Example}

We can find more SDFs of type $2$ by computer search which satisfy the hypothesis of Proposition \ref{pat2}. They are listed in Table \ref{tab:1}, which will be used to construct relative difference families later.

\begin{table}[t]
\begin{center}
\begin{tabular}{|c|l|}
\hline $(\mathbb{Z}_{12},5,20)$-SDF & $\Sigma_2=[[0, 0, 0, 6, 6],[0, 0, 0, 0, 6]]$ \\&
$\Sigma_3=\underline{2}$ $[[0,1,2,3,4],[0,1,2,4,5],[0,1,3,5,8],[0,1,4,5,8],[0,2,4,7,9]]$\\
\hline $(\mathbb{Z}_{25},6,6)$-SDF & $\Sigma_1=[[0,0,5,5,14,14]]$\\&
$\Sigma_3=\underline{2}$ $[[0,1,2,3,6,18],[0,2,8,12,15,19]]$\\
\hline $(\mathbb{Z}_{30},6,6)$-SDF & $\Sigma_1=[[0,0,6,6,16,16],[0,15,3,18,7,22]]$\\&
$\Sigma_3=\underline{2}$ $[[0,1,2,3,8,21],[0,2,5,9,13,18]]$\\
\hline $(\mathbb{Z}_{35},6,6)$-SDF & $\Sigma_1=[ [0,0,8,8,18,18]]$\\&
$\Sigma_3=\underline{2}$ $[[0,1,2,3,5,15],[0,3,7,14,23,29],[0,4,9,17,23,28]]$\\
\hline $(\mathbb{Z}_{45},6,6)$-SDF & $\Sigma_1=[[0,0,10,10,26,26]]$ \\&
$\Sigma_3=\underline{2}$ $[[0,1,3,11,17,31],[0,4,9,22,30,37],[0,1,3,7,12,25],[0,1,3,7,12,25]]$ \\
\hline $(\mathbb{Z}_5,6,12)$-SDF & $\Sigma_1=[[0,0,1,1,2,2],[0,0,2,2,4,4]]$ \\
\hline $(\mathbb{Z}_{15},6,12)$-SDF & $\Sigma_1=[[0,0,3,3,8,8],[0,0,4,4,9,9]]$ \\&
$\Sigma_3=\underline{2}$ $[[0,1,2,3,4,7],[0,1,2,4,8,10]]$ \\
\hline $(\mathbb{Z}_{35},7,6)$-SDF & $\Sigma_1=[[0, 7, 7, 17, 17, 30, 30]]$ \\&
$\Sigma_3=\underline{2}$ $[[0,1,2,3,5,21,29],[0,3,9,13,17,24,29]]$ \\
\hline $(\mathbb{Z}_{49},7,6)$-SDF & $\Sigma_1=[[ 0, 4, 4, 16, 16, 36, 36 ]]$\\&
$\Sigma_3=\underline{2}$ $[0,1,3,20,28,38,43],[0,1,3,27,31,36,42],[0,1,3,27,31,36,42]]$\\
\hline $(\mathbb{Z}_{21},7,12)$-SDF & $\Sigma_1=[[0, 5, 5, 10, 10, 17, 17],
[ 0, 3, 3, 9, 9, 17, 17]]$\\&
$\Sigma_3=\underline{2}$ $[[0,1,2,3,4,5,11],(0,1,3,7,11,13,16]]$\\
\hline
\end{tabular}
\caption{SDFs of type $2$ satisfying Proposition \ref{pat2}}\label{tab:1}
\end{center}
\end{table}

\begin{Theorem}\label{ExistenceDF2}
If there exists a $(G,k,\mu)$-SDF that satisfies the hypothesis of Proposition $\ref{pat2}$, then there exists a $(G\times \mathbb{F}_q,G\times \{0\},k,1)$-DF for any prime power $q\equiv 1 \pmod{\mu}$ and $q> Q(\mu/2,k-1)$.
\end{Theorem}

\proof Take $\mathfrak{C}=[C_1,C_2,\dots, C_n]$ as in the proof of Proposition \ref{pat2}. It has been shown that $$\Delta {\mathfrak C}=\bigcup_{g\in G} (\{g\}\times D_g),$$
where $D_g=\{1,-1\}\cdot L_g$ and $L_g$ is a multiset of size $\mu/2$ on $\mathbb{F}_q^*$. By Proposition \ref{prop:basic}, to obtain a $(G\times \mathbb{F}_q,G\times \{0\},k,1)$-DF, it suffices to prove that for each $g\in G$, $L_g$ can be taken as a system of representatives for $C_0^{\mu/2,q}$ in $\mathbb{F}_q^*$.

Note that $L_g=L_{-g}$ and each element of $L_g$ is of one of the forms: $y_{i,j_1}\pm y_{i,j_2}$, $y_{i,j}$ and $2y_{i,j}$. Therefore, by Theorem \ref{thm:cyclot bound}, for any prime power $q\equiv 1 \pmod{\mu}$ and $q> Q(\mu/2,k-1)$, it can be required that every $L_g$ is a system of representatives for $C_0^{\mu/2,q}$ in $\mathbb{F}_q^*$.
\endproof

We present the following example as an illustration of how to use Theorem \ref{ExistenceDF2} to construct DFs.

\begin{Example}\label{eg:10p}
Given any prime power $q\equiv 1 \pmod{12}$ and $q> Q(6,4)$, there exists a $(\mathbb{Z}_{10}\times \mathbb{F}_q,\mathbb{Z}_{10}\times \{0\},5,1)$-DF.
\end{Example}

\proof Take the $(\mathbb{Z}_{10},5,12)$-SDF of type $2$, $\Sigma=[A_1,A_2,\ldots,A_6]$, from Example \ref{10-5-12}, where $A_1=[0, 3, 3, 7, 7]$, $A_2=[0, 0, 0, 5, 5],$ $A_3=A_4=[0, 9, 6, 7, 8]$ and $A_5=A_6=[0, 8, 4, 6, 7].$ Consider the family $\mathfrak{C}=[C_1,C_2,\dots, C_6]$, a collection of subsets of $\mathbb{Z}_{10}\times \mathbb{F}_q$, whose first components come from $\Sigma$, where
\begin{center}
\begin{tabular}{l}
$C_1=\{(0,0), (3,y_{1,1}), (3,-y_{1,1}), (7,y_{1,2}), (7,-y_{1,2})\}$;\\
$C_2=\{(0,y_{2,1}), (0,y_{2,2}), (0,y_{2,3}), (5,y_{2,4}), (5,y_{2,5})\}$;\\
$C_3=\{(0,y_{3,1}), (9,y_{3,2}), (6,y_{3,3}), (7,y_{3,4}), (8,y_{3,5})\}$; \ \ \ $C_4= C_3\cdot\{(1,-1)\}$;\\
$C_5=\{(0,y_{4,1}), (8,y_{4,2}), (4,y_{4,3}), (6,y_{4,4}), (7,y_{4,5})\}$; \ \ \ $C_6= C_5\cdot\{(1,-1)\}$.
\end{tabular}
\end{center}
One can check that
$$\Delta [C_1,C_2,\dots, C_6]=\bigcup_{g\in \mathbb{Z}_{10}} (\{g\}\times D_g),$$
where $D_g=\{1,-1\}\cdot L_g$, $L_g=L_{-g}$ and
\begin{center}
\begin{tabular}{l}
$L_0=\{2y_{1,1},2y_{1,2},y_{2,1}-y_{2,2},y_{2,1}-y_{2,3},y_{2,3}-y_{2,2},y_{2,4}-y_{2,5}\}$;\\
$L_1=\{y_{3,2}-y_{3,1},y_{3,2}-y_{3,5},y_{3,5}-y_{3,4},y_{3,4}-y_{3,3},y_{4,2}-y_{4,5},y_{4,5}-y_{4,4}\}$;\\
$L_2=\{y_{3,1}-y_{3,5},y_{3,5}-y_{3,3},y_{3,2}-y_{3,4},y_{4,1}-y_{4,2},y_{4,2}-y_{4,4},y_{4,4}-y_{4,3}\}$;\\
$L_3=\{y_{1,1},y_{1,2},y_{3,4}-y_{3,1},y_{3,2}-y_{3,3},y_{4,1}-y_{4,5},y_{4,5}-y_{4,3}\}$;\\
$L_4=\{y_{1,2}-y_{1,1},y_{1,2}+y_{1,1},y_{4,4}-y_{4,1},y_{4,1}-y_{4,3},y_{4,3}-y_{4,2},y_{3,3}-y_{3,1}\}$;\\
$L_5=\{y_{2,4}-y_{2,1},y_{2,4}-y_{2,2},y_{2,4}-y_{2,3},y_{2,5}-y_{2,1},y_{2,5}-y_{2,2},y_{2,5}-y_{2,3}\}$.
\end{tabular}
\end{center}
By Theorem \ref{thm:cyclot bound}, for any prime power $q\equiv 1 \pmod{12}$ and $q> Q(6,4)$, we can always require that every $L_g$ is a system of representatives for $C_0^{6,q}$ in $\mathbb{F}^*_q$. Therefore, given a transversal $S$ for $\{1,-1\}$ in $C_0^{6,q}$, the family $[C\cdot \{(1,s)\} : s\in S,\ C\in \mathfrak{C}]$ is a $(\mathbb{Z}_{10}\times \mathbb{F}_q,\mathbb{Z}_{10}\times \{0\},5,1)$-DF.
\endproof

\begin{Remark}\label{newcasesI}
Let $k$ be odd and $\Sigma$ be a $(G,k,\mu)$-SDF of type $2$ satisfying the hypothesis of Proposition $\ref{pat2}$. If $\Sigma$ consists only of base blocks belonging to $\Sigma_1$, then the lower bound on $q$ in Theorem $\ref{ExistenceDF2}$ can be improved. That is to say, there exists a $(G\times \mathbb{F}_q,G\times \{0\},k,1)$-DF for any prime power $q\equiv 1 \pmod{\mu}$ and $q>Q(\mu/2,k-2)$ $($in this case every $y_{i,j}$ contributes at most $k-2$ cyclotomic conditions$)$.
\end{Remark}

\begin{Theorem}\label{DF}
Let $q$ be a prime. Then there exists a $(\mathbb{Z}_{h}\times \mathbb{F}_q,\mathbb{Z}_{h}\times \{0\},k,1)$-DF in the following cases:
\begin{center}
\begin{tabular}{|l|c|}
\hline{$(hq,h,k,1)$} & $\mbox{possible exceptions / \textbf{definite ones}}$\\
\hline $(2q,2,5,1)$: $q\equiv 1 \pmod{20}$ & \\
\hline $(10q,10,5,1)$: $q\equiv 1 \pmod{12}$ & \\
\hline $(12q,12,5,1)$: $q\equiv 1 \pmod{20}$ & \\
\hline $(hq,h,6,1)$: $h \in \{25,30,35,45\},\ q\equiv 1 \pmod{6}$ &$(25\times 7,25,6,1)$\\
\hline $(hq,h,6,1)$: $h \in \{5,15\},\ q\equiv 1 \pmod{12}$ & $ \mathbf{(5\times 13,5,6,1)}$\\
\hline $(hq,h,7,1)$: $h \in \{35,49\},\ q\equiv 1 \pmod{6}$ & $(35\times 7,35,7,1)$, $(49\times 7,49,7,1)$\\
\hline $(21q,21,7,1)$: $q\equiv 1 \pmod{12}$ &\\
\hline
\end{tabular} .
\end{center}
\end{Theorem}

\proof Begin with the SDFs of type $2$ listed in Examples \ref{eg:2,5,20-con}, \ref{10-5-12} and Table \ref{tab:1}. Then apply Theorem \ref{ExistenceDF2} and Remark \ref{newcasesI} to get the desired DFs for sufficiently large $q$. For the values of $q$ smaller than the lower bounds, we found, by computer search, all the DFs that satisfy the required cyclotomic conditions in Proposition \ref{prop:basic} except for the cases of $(h,q,k,\lambda)\in\{(10,13,5,1),(30,7,6,1),(5,13,6,1),(5$, $37,6,1),(15, 13,6,1),(25,7,6,1), (35,7,6,1), (45,7,6,1),(21,13,7,1),(35,7,7,1),(49,7,7,1)\}$. The interested reader may get a copy of these data from the authors.

For $(h,q,k,\lambda)=(10,13,5,1)$, we give an explicit construction of a $(\mathbb{Z}_{10}\times \mathbb{F}_{13},\mathbb{Z}_{10}\times \{0\},5,1)$-DF:
\begin{center}
\begin{tabular}{l}
$\{(0,0), (3,1), (3,12), (7,3), (7,10)\}$, \\
$\{(0,0), (0,1), (0,4), (5,3), (5,8)\}$, \\
$\{(0,0), (9,4), (6,10), (7,5), (8,6)\}\cdot\{(1,x)\}$, \\
$\{(0,0), (8,2), (4,7), (6,12), (7,9)\}\cdot\{(1,x)\}$.
\end{tabular}
\end{center}
where $x$ runs over $C_0^{6,13}$. For $(h,q,k,\lambda)\in\{(30,7,6,1),(5,37,6,1),(15,13,6,1),(35,7,6,1),(45,7$, $6,1),(21,13,7,1)\}$, we give an explicit construction of a $(\mathbb{Z}_{h}\times \mathbb{F}_{q},\mathbb{Z}_{h}\times \{0\},k,1)$-DF as follows:
\begin{center}
\begin{tabular}{|c|l|}
\hline $(30\times 7,30,6,1)$-DF & $\{0,11,30,111,131,171\}$, $\{0,12,27,73,148,165\}$, \\
& $\{0,1,3,25,34,128\}\times 29^i$,
$\{0,5,37,43,53,139\}\times 29^i$,  where $i=0,1$.\\\hline
$(5\times 37,5,6,1)$-DF & $\{0,4,38,95,102,138\}$, $\{0,12,29,79,109,165\}$, \\
& $\{0,1,3,11,42,123\}\times 26^i$,
$\{0,5,19,77,145,169\}\times 26^i$,  where $i=0,1$.\\\hline
$(15\times 13,15,6,1)$-DF & $\{0,1,3,8,18,89\}\times 16^i$,
$\{0,4,42,105,141,155\}\times 16^i$,  where $i=0,1,2$.\\\hline
$(35\times 7,35,6,1)$-DF & $\{0,10,30,85,130,195\}$, $\{0,1,3,54,158,167\}\times 116^i$, \\
& $\{0,4,22,93,118,201\}\times 116^i$, where $i=0,1,2$.\\\hline
$(45\times 7,45,6,1)$-DF & $\{0,1,3,11,61,244\}\times 16^i$, $\{0,4,29,135,171,278\}\times 16^i$, \\
& $\{0,5,23,78,125,164\}\times 16^i$, where $i=0,1,2$.\\\hline
$(21\times 13,21,7,1)$-DF & $\{0,1,3,9,88,116,135\}\times 16^i$,\\
& $\{0,4,37,86,127,204,211\}\times 16^i$,  where $i=0,1,2$.\\
\hline
\end{tabular}
\end{center}
For $(h,q,k,\lambda)=(5,13,6,1)$, an exhaustive search shows that no $(\mathbb{Z}_{5}\times \mathbb{F}_{13},\mathbb{Z}_{5}\times \{0\},6,1)$-DF exists. \qed

At the end of this section, we remark that Buratti \cite{b99} presented five classes of interesting SDFs. One can check that they satisfy the hypothesis of Proposition \ref{pat2}, and so they are SDFs of type 2 with respect to any odd prime power $q$.

\begin{Lemma} {\rm (implied in \cite{b99})} \label{lem:SDF-Paley}
\begin{itemize}
\item[$(1)$] Let $p$ be an odd prime power. Then the $(\mathbb{F}_p,p,p-1)$-SDF given by the single base block $\{0\}\cup \underline{2}C_0^{2,p}$ is of type $2$ $($Paley SDF of the $1$st type$)$.
\item[$(2)$] Let $p\equiv 3\pmod{4}$ be a prime power. Then the $(\mathbb{F}_p,p+1,p+1)$-SDF given by the single base block $\underline{2}(\{0\}\cup C_0^{2,p})$ is of type $2$ $($Paley SDF of the $2$nd type$)$.
\item[$(3)$] Let $p$ be an odd prime power. Set $X_1={\underline 2}(\{0\}\cup C_0^{2,p})$ and $X_2={\underline 2}(\{0\}\cup C_1^{2,p})$. Then the $(\mathbb{F}_p,p+1,2p+2)$-SDF $[X_1,X_2]$ is of type $2$ $($Paley SDF of the $3$rd type$)$.
\item[$(4)$] Given twin prime powers $p>2$ and $p+2$, the set
$$(C_0^{2,p} \times C_0^{2,p+2})\cup (C_1^{2,p}\times C_1^{2,p+2})\cup (\mathbb{F}_p\times\{0\})$$
is a $(p(p+2),\frac{p(p+2)-1}{2},\frac{p(p+2)-3}{4})$ difference set over $\mathbb{F}_p \times \mathbb{F}_{p+2}$, which is a difference family with only one base block. Let $D$ be its complement. Then ${\underline 2}D$ is a $(\mathbb{F}_p \times \mathbb{F}_{p+2},p(p+2)+1,p(p+2)+1)$-SDF of type $2$ $($twin prime power difference multiset$)$.
\item[$(5)$] Given any prime power $p$ and any integer $d$, there exists a cyclic $(\frac{p^d-1}{p-1},\frac{p^{d-1}-1}{p-1},\frac{p^{d-2}-1}{p-1})$ difference set $($Singer difference set$)$, which  is a difference family with only one base block. Let $D$ be its complement. Then ${\underline p}D$ is a $(\frac{p^d-1}{p-1},p^d,p^d(p-1))$-SDF of type $2$ $($Singer difference multiset$)$.
\end{itemize}
\end{Lemma}

\begin{Remark}\label{rek:paley type 2}
Take the Paley SDFs from Lemma $\ref{lem:SDF-Paley}(1)$ and $(2)$. Then use Theorem $\ref{ExistenceDF2}$ and Remark $\ref{newcasesI}$. We can obtain the following DFs:
\begin{itemize}
\item[1)] there exists an $(\mathbb{F}_p\times \mathbb{F}_q,\mathbb{F}_p\times \{0\},p,1)$-DF for any odd prime powers $p,q$ with $q\equiv 1 \pmod{p-1}$ and $q>Q((p-1)/2,p-2)$;
\item[2)] there exists an $(\mathbb{F}_p\times \mathbb{F}_q,\mathbb{F}_p\times \{0\},p+1,1)$-DF for any prime powers $p,q$ with $p\equiv 3 \pmod{4}$, $q\equiv 1 \pmod{p+1}$ and $q>Q((p+1)/2,p)$.
\end{itemize}
These DFs can also be found in Theorems $3.11$ and $3.12$ of $\cite{cfw}$.
\end{Remark}

\section{Strong difference families of type $4$}

In this section, we shall present a class of strong difference families of type 4, and then use them to establish new (asymptotic) existence results for DFs.

\begin{Proposition}\label{pat4}
Let $(G,+)$ be an abelian group of odd order.
Let $\Sigma$ be a $(G,k,\mu)$-SDF where $k\equiv 0,1 \pmod{4}$ and $\Sigma$ is the union of two families $\Sigma_1$ and $\Sigma_2$ $(\Sigma_2$ could be empty$)$ such that:
\begin{itemize}
\item[$1)$] $\Sigma_1$ consists of only one base block $($called the {\em distinguished base block}$)$ of the form
$$[x_{1},x_{1},-x_{1},-x_{1},\dots,x_{ k/4},x_{ k/4},-x_{ k/4},-x_{ k/4}]\hbox{ if }k\equiv 0\ ({\rm mod}\ 4),\hbox{ and}$$
$$[0,x_{1},x_{1},-x_{1},-x_{1},\dots,x_{\lfloor k/4\rfloor},x_{\lfloor k/4\rfloor},-x_{\lfloor k/4\rfloor},-x_{\lfloor k/4\rfloor}]\hbox{ if }k\equiv 1\ ({\rm mod}\ 4)$$ satisfying that $\Delta[x_1,-x_1,x_2,-x_2,\dots,x_{\lfloor k/4\rfloor},-x_{\lfloor k/4\rfloor}]$ does not contain zeros;
\item[$2)$] for any $A\in \Sigma_2$, the multiplicity of $A$ in $\Sigma_2$ is doubly even and $\Delta(A)$ does not contain zeros.
\end{itemize}
Then $\Sigma$ is of type $4$.
\end{Proposition}
\proof
Let $\Sigma=[A_1,A_2,\dots,A_n]$ be the given $(G,k,\mu)$-SDF, where $A_1=A$ is the distinguish block, $A_{2}=A_{3}=A_{4}=A_{5}$, $\dots$, $A_{n-3}=A_{n-2}=A_{n-1}=A_n$. Let $\xi$ be a primitive $4$th root of unity in $\mathbb{F}^*_q$. Consider the family $\mathfrak{C}=[C_1,C_2,\dots, C_n]$, where
\begin{itemize}
\item if $k\equiv 0 \pmod{4}$ and $A_1=[x_{1,1},x_{1,1},-x_{1,1},-x_{1,1},\dots,x_{1, k/4},x_{1, k/4},-x_{1, k/4},-x_{1, k/4}]$, then
    $$C_1=\{(x_{1,1},y_{1,1}),(x_{1,1},-y_{1,1}),(-x_{1,1},\xi y_{1,1}),(-x_{1,1},-\xi y_{1,1}),\dots, $$ $$ (x_{1, k/4},y_{1, k/4}),(x_{1, k/4},-y_{1, k/4}),(-x_{1, k/4},\xi y_{1, k/4}),(-x_{1, k/4},-\xi y_{1, k/4})\};$$
    if $k\equiv 1 \pmod{4}$ and $A_1=[0,x_{1,1},x_{1,1},-x_{1,1},-x_{1,1},\dots,x_{1,\lfloor k/4\rfloor},x_{1,\lfloor k/4\rfloor},-x_{1,\lfloor k/4\rfloor},\linebreak -x_{1,\lfloor k/4\rfloor}]$, then
    $$C_1=\{(0,0),(x_{1,1},y_{1,1}),(x_{1,1},-y_{1,1}),(-x_{1,1},\xi y_{1,1}),(-x_{1,1},-\xi y_{1,1}),\dots, $$ $$ (x_{1,\lfloor k/4\rfloor},y_{1,\lfloor k/4\rfloor}),(x_{1,\lfloor k/4\rfloor},-y_{1,\lfloor k/4\rfloor}),(-x_{1,\lfloor k/4\rfloor},\xi y_{1,\lfloor k/4\rfloor}),(-x_{1,\lfloor k/4\rfloor},-\xi y_{1,\lfloor k/4\rfloor})\};$$
\item if $A_i=A_{i+1}=A_{i+2}=A_{i+3}=[x_{i,1},x_{i,2},\dots, x_{i,k}]$ for $i\equiv 2 \pmod{4}$, then
    $$C_i=\{(x_{i,1},y_{i,1}),(x_{i,2},y_{i,2}),\dots, (x_{i,k}, y_{i,k})\},$$
    and
    $$C_{i+1}= C_i\cdot\{(1,\xi)\},\ \ C_{i+2}= C_i\cdot\{(1,-1)\},\ \ C_{i+3}= C_i\cdot\{(1,-\xi)\}.$$
\end{itemize}
Then one can check that
$$\Delta {\mathfrak C}=\bigcup_{g\in G}( \{g\}\times D_g),$$
where $D_g=\{1,-1,\xi,-\xi\}\cdot L_g$ and $L_g$ is a multiset of size $\mu/4$ on $\mathbb{F}_q^*$. Note that the zero element of $G$ must appear $\lfloor k/4\rfloor\times 4$ times in $\Delta \Sigma$, so $\mu=4\lfloor k/4\rfloor$.  It follows that $\Sigma$ is of type 4 with respect to any prime power $q\equiv 1 \pmod{4}$.
\endproof

We remark that to use Proposition \ref{pat4}, the given $(G,k,\mu)$-SDF must contain $\mu |G|/(k(k-1))\equiv 1\pmod{4}$ base blocks. Also, $\mu=4\lfloor k/4\rfloor$, and so $G$ must be a group of odd order.

\begin{Example}\label{45-5-4}
By Proposition $\ref{pat4}$, the following $(\mathbb{Z}_{45},5,4)$-SDF is of type $4$  with respect to any prime power $q\equiv 1 \pmod{4}$:
\begin{center}
\begin{tabular}{|c|c|}
\hline $\Sigma_1$ & $\Sigma_2$ \\
\hline $[[0,1,1,-1,-1]]$ & $[[0,3,7,13,30], [0,3,7,13,30], [0,3,7,13,30], [0,3,7,13,30], $\\
 & $ [0,5,14,26,34], [0,5,14,26,34], [0,5,14,26,34], [0,5,14,26,34]]$\\\hline
\end{tabular} .
\end{center}
\end{Example}

Other SDFs of type $4$ which satisfy the hypothesis of Proposition \ref{pat4} can be found in \cite{cfw}. They are listed in Table \ref{tab:2}, which will be used to construct relative difference families later.
\begin{table}[t]
\begin{center}
\begin{tabular}{|c|l|}
\hline $(\mathbb{Z}_{63},8,8)$-SDF & $\Sigma_1=[[20,20,-20,-20,29,29,-29,-29]]$ \\
& $\Sigma_2=\underline{4}$   $[[0,1,3,7,19,34,42,53],[0,1,4,6,26,36,43,51]]$\\
\hline
$(\mathbb{Z}_{81},9,8)$-SDF & $\Sigma_1=[[0,4,4,-4,-4,37,37,-37,-37]]$\\
&  $\Sigma_2=\underline{4}$ $[[0,1,4,6,17,18,38,63,72],[0,2,7,27,30,38,53,59,69]]$\\
\hline
\end{tabular}
\caption{SDFs of type $4$ satisfying Proposition \ref{pat4}}\label{tab:2}
\end{center}
\end{table}


Before stating a theorem similar to Theorem \ref{ExistenceDF2}, we present the following example as an illustration of how to apply SDFs of type 4 to obtain DFs.

\begin{Example}\label{eg:45}
Given any prime power $q\equiv 1 \pmod{4}$, there exists a $(\mathbb{Z}_{45}\times \mathbb{F}_q,\mathbb{Z}_{45}\times \{0\},5,1)$-DF.
\end{Example}

\proof
Take the $(\mathbb{Z}_{45},5,4)$-SDF, $\Sigma=[A_1,A_2,\ldots,A_9]$, from Example \ref{45-5-4}, where $A_1=[0,1,1,-1,$ $-1]$, $A_2=A_3=A_4=A_5=[0,3,7,13,30]$ and $A_6=A_7=A_8=A_9=[0,5,14,26,34]$. Let $\xi$ be a primitive $4$th root of unity in $\mathbb{F}^*_q$. Now, consider the family $\mathfrak{C}=[C_1,C_2,\dots, C_9]$ of base blocks whose first components come from $\Sigma$, where
\begin{center}
\begin{tabular}{l}
$C_1=\{(0,0),(1,y_{1,1}),(1,-y_{1,1}),(-1,y_{1,1}\xi),(-1,-y_{1,1}\xi)\}$;\\
$C_2=\{(0,y_{2,1}),(3,y_{2,2}),(7,y_{2,3}),(13,y_{2,4}),(30,y_{2,5})\}$; \\
$C_3= C_2\cdot\{(1,\xi)\};$ \ \ \ \ $C_4= C_2\cdot\{(1,-1)\};$ \ \ \ \ $C_5= C_2\cdot\{(1,-\xi)\}$;\\
$C_6=\{(0,y_{6,1}),(5,y_{6,2}),(14,y_{6,3}),(26,y_{6,4}),(34,y_{6,5})\}$;\\
$C_7= C_6\cdot\{(1,\xi)\};$ \ \ \ \ $C_8= C_6\cdot\{(1,-1)\};$ \ \ \ \ $C_9= C_6\cdot\{(1,-\xi)\}.$
\end{tabular}
\end{center}
One can check that
$$\Delta \mathfrak{C}=\bigcup_{g\in \mathbb{Z}_{45}} (\{g\}\times D_g),$$
where $D_g=\{1,-1,\xi,-\xi\}\cdot L_g$ and $|L_g|=1$ for any $g\in \mathbb{Z}_{45}$.

By Proposition \ref{prop:basic} with $\mu=d=4$, if for each $g\in \mathbb{Z}_{45}$, $L_g$ is a system of representatives for $C_0^{1,q}=\mathbb{F}_q^*$ in $\mathbb{F}_q^*$, i.e., $L_g$ consists of a nonzero element, then given a system of representatives $S$ for $C_0^{(q-1)/4,q}=\{1,-1,\xi,-\xi\}$ in $C_0^{1,q}=\mathbb{F}_q^*$, the family $[C\cdot \{(1,s)\} : s\in S,\ C\in \mathfrak{C}]$ forms a $(\mathbb{Z}_{45}\times \mathbb{F}_q,\mathbb{Z}_{45}\times \{0\},5,1)$-DF. This always can be done by taking nonzero second coordinates such that $y_{2,j_1}\neq y_{2,j_2}$ and $y_{6,j_1}\neq y_{6,j_2}$ for $1\leq j_1<j_2\leq 5$.
\endproof

We remark that a $(\mathbb{Z}_{45}\times \mathbb{F}_q,\mathbb{Z}_{45}\times \{0\},5,1)$-DF is also constructed in \cite[Lemma 2.8]{MaChang04} by using a different $(\mathbb{Z}_{45},5,4)$-SDF implicitly (see also Table \ref{tab:3}).

To prove Theorem \ref{ExistenceDF4}, we need a technical lemma. Denote by $G\setminus E$ the subgraph of a graph $G$ obtained by deleting the edges of $E$. Similarly a subgraph obtained by the deletions of the vertices of $U$ is denoted by $G-U$.

\begin{Lemma}\label{lem:DF-2-Dh2}
Let $r\geq 5$ and $\alpha$ be a permutation of the elements of $\mathbb{Z}_r$. Let $a\in \mathbb{Z}_r$. Then there exists a permutation $\pi$ of the elements of $\mathbb{Z}_r$ such that, for any $x\in \mathbb{Z}_r$ that is not fixed by $\alpha$, we have
$$\pi(\alpha(x))-\pi(x)\not=a.$$
\end{Lemma}

\proof When $a=0$, any bijection $\pi:\mathbb{Z}_r\rightarrow \mathbb{Z}_{r}$ can lead to the desired conclusion. Assume that $a\neq 0$.

We first give an equivalent description of this lemma in the language of graphs. We define a directed graph $\overrightarrow{H}$, whose vertices are taken from $\mathbb{Z}_r\setminus \{x\in \mathbb{Z}_r:x=\alpha(x)\}$. For any two distinct vertices $x$ and $y$, $(x,y)$ is a directed edge of $\overrightarrow{H}$ if and only if $y=\alpha(x)$. We have that $\overrightarrow{H}$ is the union of disjoint directed cycles of lengths $l_1,l_2,\dots,l_s$, where each $l_i\geq 2$ and $\sum_{i=1}^s l_i\leq r$. Now regard the mapping $\pi$ as an embedding of $\overrightarrow{H}$ in the complete directed graph $\overrightarrow{K_r}$ defined on $\mathbb{Z}_{r}$. Consider the directed graph $\overrightarrow{C}$ whose edge set is given by the pairs $\{(x,x+a): x \in \mathbb{Z}_r\}$. $\overrightarrow{C}$ is the union of directed cycles whose length $l$ is the order of $a$ in $\mathbb{Z}_r$.
Now we prove that for any $r\geq 5$, $\overrightarrow{H}$ can be embedded in $\overrightarrow{K_r}\setminus \overrightarrow{C}$.

CASE 1: $a=r/2$.
In the case no length $l_i=r/2+1$ we can embed $\overrightarrow{H}$ in $\overrightarrow{K_r}\setminus \overrightarrow{C}$ as the union of the cycles $(0,\dots,l_1-1),(l_1,\dots,l_1+l_2-1),\dots, ((\sum_{i=1}^{s-1} l_i),\dots,(\sum_{i=1}^{s} l_i)-1 )$.
Otherwise since $\sum_{i=1}^s l_i\leq r$, there is exactly one length $l_i$ that is equal to $r/2+1$ and we can assume it is $l_1$. In this case we can embed $\overrightarrow{H}$ as the union of the cycles $(0,1,\dots,l_1-2,l_1), (l_1-1, l_1+1,\dots,l_1+l_2-1), (l_1+l_2,l_1+l_2+1,\dots,l_1+l_1+l_3-1),\dots, ((\sum_{i=1}^{s-1} l_i),\dots,(\sum_{i=1}^{s} l_i)-1 ).$
Since here $r$ is even, $r\geq 6$ and so $a\neq 2$, which makes the directed edges $(l_1-2,l_1)$ and $(l_1-1, l_1+1)$ admissible.

CASE 2: $a\not=0,r/2$. In this case $\overrightarrow{C}$ is the union of oriented cycles of length $l$.
Assume that $l_1=l_2=\dots=l_{s'}=2$ and $l_{i}\neq 2$ for $i>s'$.
Since $r\geq 5$, the unoriented graph $K_r\setminus \overrightarrow{C}$ admits a $1$-factor for even values of $r$ and an almost $1$-factor for odd values. It follows that we can embed the union $\overrightarrow{H_1}$ of $s'$ cycles of length $2$ in $\overrightarrow{K_r}\setminus \overrightarrow{C}$.
We note that $\overrightarrow{C}- \overrightarrow{H_1}$ is the union of oriented paths $\overrightarrow{P}_1,\overrightarrow{P}_2,\dots ,\overrightarrow{P}_p$ and cycles $\overrightarrow{C}_1,\overrightarrow{C}_2,\dots, \overrightarrow{C}_c$.
We denote the vertices of $\overrightarrow{C}- \overrightarrow{H_1}$
by $v_1,v_2,\dots,v_{r-2s'}$ in such a way that, if the path $\overrightarrow{P}_i$ has length $b_i$, it is of the form $[v_{h_i+b_i},v_{h_i+b_i-1},v_{h_i+b_i-2},\dots,v_{h_i}]$ for some integer $h_i$ and the cycle $\overrightarrow{C}_i$ is of the form $(v_{u_i+l},v_{u_i+l-1},\dots,v_{u_i+1})$ for some integer $u_i$.
Now, since $a\not=r/2$, we can embed $\overrightarrow{H}_2=\overrightarrow{H}- \overrightarrow{H}_1$ in $(\overrightarrow{K_r}- \overrightarrow{H_1})\setminus \overrightarrow{C}$ as the union of the cycles
$$(v_1,\dots,v_{l_{s'+1}}),(v_{l_{s'+1}+1},\dots,v_{l_{s'+1}+l_{s'+2}}),\dots, $$ $$(v_{l_{s'+1}+l_{s'+2}+\dots+l_{s-1}+1},\dots,v_{l_{s'+1}+l_{s'+2}+\dots+l_{s}}).$$
This completes the proof.
\endproof

\begin{Theorem}\label{ExistenceDF4}
If there exists a $(G,k,\mu)$-SDF that satisfies the hypothesis of Proposition $\ref{pat4}$, whose distinguished base block is denoted by $A$, then for any prime power $q\equiv 1 \pmod{\mu}$ and $q> Q(\mu/4,k-1)$, there exists a $(G\times \mathbb{F}_q,G\times \{0\},k,1)$-DF in the following cases
\begin{itemize}
\item[$1)$] $k\equiv 0 \pmod{4}$;
\item[$2)$] $k=5$;
\item[$3)$] $k\in \{9,13,17\}$ and $x_{1}\not=\pm 2x_{2}$ for any nonzero $x_1,x_2\in A$ $(x_1$ could be $x_2)$;
\item[$4)$] $k\equiv 1 \pmod{4}$, $k\geq 21$ and $3x\not=0$ for any nonzero $x\in A$.
\end{itemize}
\end{Theorem}

\proof Take $\mathfrak{C}=[C_1,C_2,\dots, C_n]$ as in the proof of Proposition \ref{pat4}. It has been shown that $$\Delta {\mathfrak C}=\bigcup_{g\in G}( \{g\}\times D_g),$$
where $D_g=\{1,-1,\xi,-\xi\}\cdot L_g$ and $L_g$ is a multiset of size $\mu/4$ on $\mathbb{F}_q^*$. It is readily checked that $L_g=L_{-g}$ and each element of $L_g$ is one of the following forms: $y_{1,j}$, $2y_{1,j}$, $(1-\xi)y_{1,j}$, $y_{1,j_1}\pm \xi y_{1,j_2}$, $y_{1,j_1}\pm y_{1,j_2}$ and $y_{i,j_1}-y_{i,j_2}$ for $i\geq 2$. Note that $2y_{1,j}\in L_0,$ $(1-\xi)y_{1,j}\in L_{\pm 2x_{1,j}}$ and when $k\equiv 1\pmod{4}$, $y_{1,j}\in L_{\pm x_{1,j}}$. By Proposition \ref{prop:basic}, to obtain a $(G\times \mathbb{F}_q,G\times \{0\},k,1)$-DF, it suffices to prove that for each $g\in G$, $L_g$ can be taken as a system of representatives for $C_0^{\mu/4,q}$ in $\mathbb{F}_q^*$. By Theorem \ref{thm:cyclot bound}, this can be done for any prime power $q\equiv 1 \pmod{\mu}$ and $q> Q(\mu/4,k-1)$ if no $L_g$ contains a $2$-subset of the form $\{y_{1,j_1},(1-\xi)y_{1,j_2}\}$.

When $k\equiv 0\pmod{4}$, $y_{1,j_1}$ cannot occur in any $L_g$, so no $L_g$ contains a $2$-subset of the form $\{y_{1,j_1},(1-\xi)y_{1,j_2}\}$.

When $k=5$, $|L_g|=1$ for any $g\in G$, so no $L_g$ contains a $2$-subset of the form $\{y_{1,j_1},(1-\xi)y_{1,j_2}\}$.

When $k\equiv 1\pmod{4}$, $k\geq 9$ and some $L_g$ contains a $2$-subset of the form $\{y_{1,j_1},(1-\xi)y_{1,j_2}\}$, since $y_{1,j_1}\in L_{\pm x_{1,j_1}}$ and $(1-\xi)y_{1,j_2}\in L_{\pm 2x_{1,j_2}}$, we have $g=\pm x_{1,j_1}=\pm 2x_{1,j_2}$. Thus if $x_{1,j_1}\not=\pm 2x_{1,j_2}$ for any nonzero $x_{1,j_1},x_{1,j_2}\in A_1$ ($j_1$ could be $j_2$), then no $L_g$ contains a $2$-subset of the form $\{y_{1,j_1},(1-\xi)y_{1,j_2}\}$.

Assume that $k\equiv 1\pmod{4}$ and $k\geq 21$. We shall show that even if some $L_g$ contains a $2$-subset of the form $\{y_{1,j_1},(1-\xi)y_{1,j_2}\}$, we can still require every $L_g$, $g\in G$, is a system of representatives for $C_0^{\mu/4,q}$ in $\mathbb{F}^*_q$, provided for any nonzero $x\in A$, $3x\not=0$.

Let $P$ be the set of ordered pairs $(y_{1,j_1},y_{1,j_2})$ for all possible $j_1$ and $j_2$ such that $(y_{1,j_1},(1-\xi)y_{1,j_2})$ is a $2$-subset of some $L_g$. We note that, if $(y_{1,j_1},y_{1,j_2})$ and $(y_{1,j_1},y_{1,j_3})\in P$, then $x_{1,j_1}=\pm 2x_{1,j_2}=\pm 2x_{1,j_3}$ and hence $2(x_{1,j_2}\pm x_{1,j_3})=0$ but, since $\Delta[x_{1,1},-x_{1,1},\dots,x_{1,\lfloor k/4\rfloor},-x_{1,\lfloor k/4\rfloor}]$ does not contain involutions ($G$ is of odd order) and zeros this is possible only if $j_2=j_3$. It follows that the multiset $P_1:=[y_{1,j_1} : (y_{1,j_1},y_{1,j_2})\in P]$ is a simple set and hence we can define a map $\alpha$ from $P_1$ to $P_2:=[y_{1,j_2} \mid (y_{1,j_1},y_{1,j_2})\in P]$ such that $\alpha(y_{1,j_1})=y_{1,j_2}$.
We also note that, if $(y_{1,j_2},y_{1,j_1})$ and $(y_{1,j_3},y_{1,j_1})\in P$, then $x_{1,j_2}=\pm 2x_{1,j_1}=\pm x_{1,j_3}$ but this is possible only if $j_2=j_3$. It follows that the multiset $P_2$ is also a simple set and hence $\alpha$ is an injective map. Moreover, if $(y_{1,j_1},y_{1,j_2})\in P$ and $j_1=j_2$, then $x_{1,j_1}=\pm 2x_{1,j_1}$ and hence $3x_{1,j}=0$. It is impossible since for any nonzero $x\in A$, $3x\not=0$ and hence $\alpha$ has no fixed point in $P_1$.

Let $r=\mu/4=\lfloor k/4\rfloor\geq 5$. Let $1-\xi\in C_{r-\alpha}^{r,q}$. Since $L_0=\{2y_{1,1},2y_{1,2},\ldots,2y_{1,r}\}$, which is a system of representatives for $C_0^{r,q}$ in $\mathbb{F}^*_q$, we can regard $P_1$ and $P_2$ as subsets of $\mathbb{Z}_r$ and we can extend the map $\alpha$ to a permutation of the elements of $\mathbb{Z}_r$. By Lemma \ref{lem:DF-2-Dh2}, there exists a bijection $\pi:\{y_{1,1},y_{1,2},\ldots,y_{1,r}\}\rightarrow \mathbb{Z}_{r}$ such that $\pi(y_{1,j_1})\not\equiv \pi(y_{1,j_2})+\alpha\pmod{r}$ for any $(y_{1,j_1},y_{1,j_2})\in P$. Thus we can assign to each $y_{1,j}$ the cyclotomic class given by the map $\pi$.
\endproof

\begin{Remark}\label{newcasesII}
Let $\Sigma$ be a $(G,k,\mu)$-SDF of type $4$ satisfying the hypothesis of Proposition $\ref{pat4}$. If $\Sigma=\Sigma_1$, i.e., $\Sigma$ consists only of the distinguished base block, then every $y_{1,j}$ contributes at most $k-3$ cyclotomic conditions when $k\equiv 0 \pmod{4}$ $($resp. $k-4$ when $k\equiv 1 \pmod{4})$. Therefore, the lower bound on $q$ in Theorem $\ref{ExistenceDF4}$ can be improved, that is to say, $q>Q(\mu/4,k-3)$ when $k\equiv 0 \pmod{4}$ and $q>Q(\mu/4,k-4)$ when $k\equiv 1 \pmod{4}$.
\end{Remark}

\begin{Theorem}\label{DF1}
Let $q$ be a prime. Then there exists a $(\mathbb{Z}_{h}\times \mathbb{F}_q,\mathbb{Z}_{h}\times \{0\},k,1)$-DF in the following cases:
\begin{center}
\begin{tabular}{|l|c|}
\hline{$(hq,h,k,1)$} & $\mbox{ possible exceptions }$\\
\hline $(63q,63,8,1)$-DF: $q\equiv 1 \pmod{8}$ & $(63\times 17,63,8,1)$\\
\hline $(81q,81,9,1)$-DF: $q\equiv 1 \pmod{8}$ & $(81\times 17,81,9,1)$, $(81\times 41,81,9,1)$\\
\hline
\end{tabular}
\end{center}
\end{Theorem}

\proof Start from the SDFs of type $4$ listed in Table \ref{tab:2}. Then apply Theorem \ref{ExistenceDF4} to get the desired DFs for sufficiently large $q$. For the values of $q$ smaller than the lower bounds, we found, by computer search, all the DFs that satisfy the required cyclotomic conditions in Proposition \ref{prop:basic} except for the cases of $(h,q,k,\lambda)\in\{(63,17,8,1),(63,41,8,1),(81,17,9,1),(81,41,9,1)\}$. The interested reader may get a copy of these data from the authors. For $(h,q,k,\lambda)=(63,41,8,1)$, we give here an explicit construction for a $(\mathbb{Z}_{63}\times \mathbb{F}_{41},\mathbb{Z}_{63}\times \{0\},8,1)$-DF:
\begin{center}
\begin{tabular}{l}
$\{(20,0),(20,1),(-20,7),(-20,35),(29,5),(29,37),(-29,18),(-29,24)\}\cdot\{(1,x)\};$\\
$\{(0,0),(1,1),(3,7),(7,4),(19,2),(34,3),(42,6),(53,27)\}\cdot\{(1,y)\};$\\
$\{(0,0),(1,3),(4,2),(6,1),(26,8),(36,29),(43,36),(51,15)\}\cdot\{(1,y)\},$\\
\end{tabular}
\end{center}
where $x$ runs over $C_0^{8,41}$ and $y$ runs over $C_0^{2,41}$. \qed



\begin{Remark}\label{rmk:paley type 4}
The Paley $(p,p,p-1)$-SDFs of the $1$st type in Lemma $\ref{lem:SDF-Paley}(1)$ satisfy the hypothesis of Proposition $\ref{pat4}$ when $p\equiv 1 \pmod{4}$ since $-1\in C_0^{2,p}$. Then by Remark $\ref{newcasesII}$, we have that there exists an $(\mathbb{F}_p\times \mathbb{F}_q,\mathbb{F}_p\times \{0\},p,1)$-DF for any prime powers $p$ and $q$ with $p\equiv 1,5\pmod{12}$, $p\neq 17$, $q\equiv 1 \pmod{p-1}$ and $q>Q((p-1)/4,p-4)$ $($see also Theorem $3.8(1)$ of $\cite{cfw})$.
\end{Remark}

\begin{Remark}\label{rmk:compare}
Comparing Remark $\ref{rek:paley type 2}(1)$ with Remark $\ref{rmk:paley type 4}$, we can see that the lower bound on $q$ in Remark $\ref{rmk:paley type 4}$ is $Q((p-1)/4,p-4)$, while it is $Q((p-1)/2,p-2)$ in Remark $\ref{rek:paley type 2}$. Hence in some circumstances, SDFs of type $4$ can produce better asymptotic existence results on DFs than SDFs of type $2$.
\end{Remark}

\section{Applications}

\subsection{A class of group divisible designs}

Difference families can be used to construct group divisible designs. Let $K$ be a set of positive integers. A {\em group divisible design} (GDD) $K$-GDD is a triple ($X, {\cal G},{\cal A}$)
satisfying the following properties: ($1$) $\cal G$ is a partition of a finite set $X$ into subsets (called {\em groups}); ($2$) $\cal A$ is a set of subsets of $X$ (called {\em blocks}), whose cardinalities are from $K$, such that every $2$-subset of $X$ is either contained in exactly one block or in exactly one group, but not in both. If $\cal G$ contains $u_i$ groups of size $g_i$ for $1\leq
i\leq r$, then $g_1^{u_1}g_2^{u_2}\cdots g_r^{u_r}$ is called the {\em type} of the GDD. The notation $k$-GDD is used when $K=\{k\}$.

\begin{Lemma}\label{6-GDD}{\rm \cite[Table 3.18]{ag}}
There exists a $6$-GDD of type $30^u$ for $u\in\{6,16,21,26,31,36,41,51,61$, $66,71,76, 78,81,86,90,91,96\}.$
\end{Lemma}

An {\em automorphism group} of a GDD $(X, {\cal G},{\cal A})$ is a permutation group on $X$ leaving ${\cal G}$ and ${\cal A}$ invariant, respectively.

Given a $(g,n,k,\lambda)$-DF, $\mathfrak{B}=[B_1,B_2,\dots,B_r]$, over an abelian group $(G,+)$ of order $g$ relative to a subgroup $N$ of order $n$, we define Dev$(B_1,B_2,\dots,B_r)$ to be the multiset
$\bigcup_{i=1}^r \{\{x+g:x\in B_i\}:g\in G\}.$
Let ${\cal N}$ be the set of cosets of $N$ in $G$. Then $(G,{\cal N}$, Dev$(B_1,B_2,\dots,B_r))$ is a $k$-GDD of type $n^{g/n}$ and its automorphism group contains a subgroup that is isomorphic to $G$.

\begin{Lemma}\label{6-GDD-1}
There exists a $6$-GDD of type $30^q$ admitting $\mathbb{Z}_{30}\times \mathbb{Z}_q$ as an automorphism group for any prime $q\equiv 1\pmod{6}$.
\end{Lemma}

\proof By Theorem \ref{DF}, there exists a $(\mathbb{Z}_{30}\times \mathbb{Z}_q,\mathbb{Z}_{30}\times \{0\},6,1)$-DF for any prime $q\equiv 1\pmod{6}$. Develop its base blocks under $\mathbb{Z}_{30}\times \mathbb{Z}_q$ to obtain a $6$-GDD of type $30^q$. \qed

\begin{Lemma}\label{6-GDD-2}
There exists a $6$-GDD of type $30^u$ for $u\in \{42,48,84,85\}$.
\end{Lemma}

\proof It is known that there exists a $7$-GDD of type $g^7$ for $g\in\{7,13\}$, which is equivalent to $5$ mutually orthogonal Latin squares of order $g$ (see Table 3.87 in \cite{ACD}). Start from a $7$-GDD of type $7^7$, and delete 7 elements from the same block to obtain a $\{6,7\}$-GDD of type $6^7$, or delete one element to obtain a $\{6,7\}$-GDD of type $7^66^1$. Start from a $7$-GDD of type $13^7$, and delete 7 elements in one group to obtain a $\{6,7\}$-GDD of type $13^66^1$, or delete 6 elements in one group to obtain a $\{6,7\}$-GDD of type $13^67^1$.

Begin with the four $\{6,7\}$-GDDs and give weight 30 to each element. Then apply Wilson's Fundamental Construction (see Theorem 2.5 in \cite{gm}) to obtain $6$-GDDs of types $180^7$, $210^6 180^1$, $390^6180^1$ and $390^6 210^1$, where the needed $6$-GDD of type $30^6$ is from Lemma \ref{6-GDD}, and the needed $6$-GDD of type $30^7$ is from Lemma \ref{6-GDD-1}.

Finally, fill in groups with $6$-GDDs of types $30^6$, $30^7$ and $30^{13}$ (from Lemma \ref{6-GDD-1}). Then we get a $6$-GDD of type $30^u$ for $u\in \{42,48,84,85\}$. \endproof

\begin{Lemma} \label{6-GDD-3}
There exists a $6$-GDD of type $30^q$ for $q\in \{25,49\}$.
\end{Lemma}

\proof Since $q$ is a prime power, we can construct a $(\mathbb{Z}_{30}\times \mathbb{F}_q,\mathbb{Z}_{30}\times \{0\},6,1)$-DF such that it satisfies the required cyclotomic conditions in Proposition \ref{prop:basic}, and then develop its base blocks under $\mathbb{Z}_{30}\times \mathbb{F}_q$ to obtain a $6$-GDD of type $30^q$. The $(\mathbb{Z}_{30},6,6)$-SDF from Table \ref{tab:1} is of type $2$ with respect to $q\in \{25,49\}$. For $q=25$, take $x^2-x+2$ to be a primitive polynomial of degree $2$ over $\mathbb{F}_5$ and $\omega$ to be a primitive root in $\mathbb{F}_{25}$. Let $\mathfrak C=\{C_1,C_2,\ldots,C_6\}$, where
\begin{center}
\begin{tabular}{l}
$C_1=\{(0,1),(0,-1),(6,\omega^{13}),(6,-\omega^{13}),(16,\omega^{11}),(16,-\omega^{11})\}$;\\
$C_2=\{(0,1),(15,-1),(3,\omega^{16}),(18,-\omega^{16}),(7,\omega^{5}),(22,-\omega^{5})\}$;\\
$C_3=\{(0,0),(1,1),(2,\omega^{20}),(3,\omega^{22}),(8,\omega^{18}),(21,\omega^{12})\}$;\ \ \ $C_4= C_3\cdot\{(1,-1)\} $;\\
$C_5=\{(0,0),(2,1),(5,\omega^{6}),(9,\omega^{19}),(13,\omega^{14}),(18,\omega^{4})\}$;\ \ \ \ $C_6= C_5\cdot\{(1,-1)\}$.\\
\end{tabular}
\end{center}
For $q=49$, take $x^2-x+3$ to be a primitive polynomial of degree $2$ over $\mathbb{F}_7$ and $\omega$ to be a primitive root in $\mathbb{F}_{49}$. Let $\mathfrak C=\{C_1,C_2,\ldots,C_6\}$, where
\begin{center}
\begin{tabular}{l}
$C_1=\{(0,1),(0,-1),(6,\omega^{40}),(6,-\omega^{40}),(16,\omega^{26}),(16,-\omega^{26})\}$;\\
$C_2=\{(0,1),(15,-1),(3,\omega^{13}),(18,-\omega^{13}),(7,\omega^{5}),(22,-\omega^{5})\}$;\\
$C_3=\{(0,0),(1,1),(2,\omega^{8}),(3,\omega^{24}),(8,\omega^{21}),(21,\omega^{10})\}$;\ \ \ $C_4= C_3\cdot\{(1,-1)\} $;\\
$C_5=\{(0,0),(2,1),(5,\omega^{34}),(9,\omega^{32}),(13,\omega^{2}),(18,\omega^{8})\}$;\ \ \ $C_6= C_5\cdot\{(1,-1)\}$.\\
\end{tabular}
\end{center}
It is readily checked that $D_g=\{1,-1\}\cdot L_g$ for $g=0,1,\ldots,29$, where $L_g$ is a system of representatives for $C_0^{3,q}$ in $\mathbb{F}_{q}^*$. Let $S$ be a system of representatives for $C_0^{(q-1)/2,q}=\{1,-1\}$ in $C_0^{3,q}$. Then by Proposition \ref{prop:basic}, $[C\cdot \{(1,s)\} : s\in S,\ C\in \mathfrak{C}]$ forms a $(\mathbb{Z}_{30}\times \mathbb{F}_q,\mathbb{Z}_{30}\times \{0\},6,1)$-DF for $q\in \{25,49\}$. \qed

Combining results of Lemmas \ref{6-GDD}, \ref{6-GDD-1}, \ref{6-GDD-2} and \ref{6-GDD-3}, we have the following theorem.

\begin{Theorem}
There exists a $6$-GDD of type $30^u$ for $u\in \{6,16,21,25,26,36,41,42,48,49,51,66,71$, $76,78, 81,84,85,86,90,91,96\}\cup\{q:q\equiv 1\ {({\rm mod}}\ 6)\ {\rm is}\ {\rm a}\ {\rm prime}\}.$
\end{Theorem}

\subsection{A class of $r$-rotational balanced incomplete block designs}

If each group of a $k$-GDD ($X, {\cal G},{\cal A}$) contains only one element and $|X|=v$, then the $k$-GDD is referred to as a {\em balanced incomplete block design}, denoted by a $(v,k,1)$-BIBD. A $(v,k,1)$-BIBD is said to be {\em $r$-rotational}, if it admits an automorphism consisting of one fixed point and $r$ cycles of length $(v-1)/r$.

\begin{Lemma} \label{lemma-k-rotational}
There exists an $r$-rotational $(31,6,1)$-BIBD for $r\in\{6,10\}$.
\end{Lemma}

\proof There is only one $(31,6,1)$-BIBD up to isomorphism $($see Table $6.5$ in $\cite{ko})$, which is a projective plane of order $5$. Check its full automorphism group and one can see that it admits an automorphism consisting of one fixed point and $r$ cycles of length $30/r$ for each $r\in\{6,10\}$.
In order to facilitate the readers, we here list all blocks of a $(31,6,1)$-BIBD over $\mathbb{Z}_{31}$: $\{i, 1+i, 3+i, 8+i, 12+i, 18+i\}$, $i\in \mathbb{Z}_{31}$. It is readily checked that the permutations in $\mathbb{Z}_{31}$
$$(1\ 13\ 6)(2\ 23\ 29)(3\ 25\ 20)(4\ 5\ 30)(7\ 26\ 28)(8\ 16\ 19)(9\ 17\ 10)(11\ 24\ 15)(12\ 21\ 27)(14\ 22\ 18)$$
and
$$(0\ 18\ 29\ 7\ 22)(1\ 25\ 26\ 16\ 19)(2\ 11\ 28\ 9\ 24)(3\ 4\ 27\ 12\ 20)(5\ 6\ 8\ 17\ 13)(10\ 15\ 14\ 21\ 30)$$ are its automorphisms of orders $3$ and $5$, respectively. \qed

\begin{Theorem} \label{BIBD-k-rotational}
Let $r\in\{6,10\}$. There exists an $r$-rotational $(30q+1,6,1)$-BIBD for any prime $q\equiv 1\pmod{6}$.
\end{Theorem}

\proof By Lemma \ref{6-GDD-1}, there exists a $6$-GDD of type $30^q$ admitting $\mathbb{Z}_{30}\times \mathbb{Z}_q$ as an automorphism group for any prime $q\equiv 1\pmod{6}$. This GDD is defined on $\mathbb{Z}_{30}\times \mathbb{Z}_q$ with groups $\mathbb{Z}_{30}\times \{i\}$, $i\in \mathbb{Z}_{q}$. Denote by $\cal A$ the set of its blocks.

For each $i\in \mathbb{Z}_{q}$, by Lemma \ref{lemma-k-rotational}, we can construct an $r$-rotational $(31,6,1)$-BIBD for $r\in\{6,10\}$ on $(\mathbb{Z}_{30}\times \{i\})\cup\{\infty\}$ such that the map $(x,i)\mapsto (x+r,i)$ is an automorphism, where $\infty$ is a fixed point. Denote by ${\cal B}_i$ the set of its blocks. Without loss of generality, we require
$${\cal B}_{i+1}={\cal B}_i+1:=\{\{(x,i+1):(x,i)\in B\}\cup\{\infty:\infty\in B\}:B\in{\cal B}_{i}\}.$$
Then ${\cal A}\cup(\bigcup_{i=0}^{q-1} {\cal B}_i)$ forms the set of blocks of a $(30q+1,6,1)$-BIBD.

Clearly, the maps $\alpha:(x,i)\mapsto (x+r,i)$ and $\beta:(x,i)\mapsto (x,i+1)$ ($\infty$ is fixed) are both automorphisms of the resulting BIBD. Since $30/r$ and $q$ are coprime, $<\alpha,\beta>$ is an automorphism group of order $30q/r$, which contains an automorphism consisting of $r$ cycles of length $30q/r$. Therefore, the resulting BIBD is $r$-rotational. \qed

\subsection{Optical orthogonal codes}

Finally we apply Theorems \ref{DF} and \ref{DF1} to construct optimal optical orthogonal codes.
A $(v,k,1)$-optical orthogonal code (OOC) is defined as a set of $k$-subsets (called {\em codewords}) of $\mathbb{Z}_v$ whose list of differences does not contain repeated elements. It is {\em optimal} if the size of the set of missing differences is less than or equal to $k(k-1)$.

A cyclic $(gv,g,k,1)$-DF can be seen as a $(gv,k,1)$-OOC whose set of missing differences is $\{0,v,2v,\ldots,(g-1)v\}$. Furthermore, one can construct a $(g,k,1)$-OOC on the set of missing differences to produce a new $(gv,k,1)$-OOC.

\begin{Lemma} \label{con:OOC} {\rm (cf. \cite[Construction 4.1]{Yin98})}
If there exist a cyclic $(gv,g,k,1)$-DF and an optimal $(g,k,1)$-OOC, then there exists an optimal $(gv,k,1)$-OOC.
\end{Lemma}

\begin{Theorem}
\begin{itemize}
\item[$(1)$] There exist an optimal $(2q,5,1)$-OOC and an optimal $(12q,5,1)$-OOC for any prime $q\equiv 1\ ({\rm mod}\ 20)$.
\item[$(2)$] There exists an optimal $(gq,k,1)$-OOC where $(g,k)\in\{(10,5),(5,6),(15,6),(21,7)\}$ for any prime $q\equiv 1\ ({\rm mod}\ 12)$ except for $(g,q,k)=(5,13,6)$.
\item[$(3)$]There exists an optimal $(gq,k,1)$-OOC where $(g,k)\in\{(25,6),$ $(30,6),$ $(35,6),$ $(45,6),$ $(35,7),$ $(49,7)\}$ for any prime $q\equiv 1\ ({\rm mod}\ 6)$ except for $(g,q,k)\in\{(25,7,6),(35,7,7),(49,7,7)\}$.
\item[$(4)$] There exists an optimal $(63q,8,1)$-OOC for any prime $q\equiv 1\ ({\rm mod}\ 8)$ and $q>17$.
\end{itemize}
\end{Theorem}

\proof Take cyclic $(hq,h,k,1)$-DFs from Theorems \ref{DF} and \ref{DF1} (note that when $h$ and $q$ are coprime, $\mathbb{Z}_{h}\times \mathbb{Z}_q$ is isomorphic to $\mathbb{Z}_{hq}$). Then apply Lemma \ref{con:OOC} with the following optimal $(h,k,1)$-OOC:
\begin{center}
\begin{tabular}{|c|c|}\hline
$(g,k,1)$ & codewords \\
\hline $(2,5,1)$, $(10,5,1)$, $(12,5,1)$ & $\emptyset$ \\
\hline $(5,6,1)$, $(15,6,1)$, $(25,6,1)$, $(30,6,1)$ & $\emptyset$ \\
\hline $(21,7,1)$, $(35,7,1)$ & $\emptyset$ \\
\hline $(35,6,1)$, $(45,6,1)$ & $\{0,1,3,7,12,20\}$ \\
\hline $(49,7,1)$  & $\{0,1,3,7,27,35,40\}$ \\
\hline $(63,8,1)$  & $\{0,1,3,7,15,20,31,41\}$ \\\hline
\end{tabular} .
\end{center}
Note that each codeword in an $(h,k,1)$-OOC contributes $k(k-1)$ different nonzero differences, so when $h\leq k(k-1)$, an optimal $(h,k,1)$-OOC is without any codeword. \qed

We remark that there is no $(81,9,1)$-OOC with one codeword by exhaustive search. Hence we cannot employ cyclic $(81q,81,9,1)$-DFs from Theorem \ref{DF1} to get optimal OOCs. For more information on combinatorics aspects of OOCs, we refer the reader to \cite{b4,cj,cm,ChangYin04,fmi,gy,MaChang04,MaChang05,Yin98} for example.

\section{Concluding remarks}

By introducing strong difference families of special types explicitly, we construct new relative difference families in this paper from the point of view of both asymptotic existences and concrete examples. Strong difference families have been widely used implicitly in the literature. We report some of them in Table \ref{tab:3} (note that here some SDFs of type $d$, $d\in \{2,4\}$, do not satisfy Propositions \ref{pat2} and \ref{pat4}).

Many infinite classes of relative difference families were obtained by using SDFs (in some cases implicitly) which do not satisfy Definition 2, see for instance \cite{b99,bp,ChangYin04,gy,MaChang04,MaChang05,m}. A further direction is to generalize Definition 2 to cover some of those constructions.

\begin{table}
\begin{center}
\begin{tabular}{|c|c|l|}
\hline Type $2$ & Source & Base blocks \\
\hline $(\mathbb{Z}_9,4,8)$-SDF & \cite[Lemma 2.11]{fmi}& $
[0,0,5,5],[0,0,3,3],\underline{4}[0,1,3,8]$.\\
\hline $(\mathbb{Z}_{12},4,4)$-SDF & \cite[Theorem 6.2]{gy} & $[0,0,5,5],
[0,2,6,8],
[0,1,3,4],
[0,8,9,11].$\\

\hline $(\mathbb{Z}_4,5,20)$-SDF & \cite[Theorem 3.3]{byw} & $[0,0,1,1,2],[0,0,0,1,2],[3,3,2,2,1]$,\\&& $[3,3,3,2,1].$\\
\hline $(\mathbb{Z}_{20},5,12)$-SDF & \cite[Theorem 2.4]{yyl} & $\underline{6}[0,1,3,9,14],\underline{2}[0,1,4,5,8],\underline{2}[0,1,5,13,18]$,\\ && $[0,2,2,18,18],[0,0,0,10,10].$ \\


\hline $(\mathbb{Z}_{45},9,8)$-SDF & \cite[Lemma 2.10]{cfw} & $[0, 2, 2, 15, 15, 23, 23, 33, 33]$,\\ && $\underline{2}[0, 1, 4, 5, 6, 7, 13, 22, 33],$\\ & & $ \underline{2}[0, 2, 5, 11, 21, 25, 28, 36, 40].$\\
\hline $(\mathbb{Z}_{125},6,6)$-SDF & \cite[Lemma 13]{cfw2} & $[0, 0, 19, 19, 71, 71],\underline{2}[0, 10, 28, 51, 78, 97]$,\\
&& $\underline{2}[0, 3, 62, 75, 86, 110],$ $ \underline{2}[0, 5, 12, 58, 70, 112]$, \\
&& $\underline{2}[0, 7, 27, 44, 70, 96], \underline{2}[0, 1, 42, 93, 85, 45],$\\
&& $\underline{2}[0, 1, 100, 104, 109, 88], \underline{2}[0, 1, 90, 81, 21, 32]$,\\
&& $\underline{2}[0, 3, 16, 40, 46, 50], $ $\underline{2}[0, 2, 7, 29, 35, 68]$, \\
&& $\underline{2}[0, 2, 8, 57, 102, 116], $ $\underline{2}[0, 2, 22, 32, 36, 96],$\\
&& $\underline{2}[0, 8, 23, 38, 72, 86].$\\\hline&&\\
\hline Type $3$ &&\\
\hline $(\mathbb{Z}_4,4,12)$-SDF & \cite[Case 1 in Section 2]{b4} & $\underline{3}[0,0,1,3],[0,2,2,2].$\\
\hline $(\mathbb{Z}_4,4,18)$-SDF & \cite[Theorem 3.2]{b4} & $\underline{3}[0,1,2,3],\underline{2}[0,0,0,1],[0,0,0,2].$\\\hline&&\\

\hline Type $4$ &&\\
\hline $(\mathbb{Z}_6,4,8)$-SDF & \cite[Case 2 in Theorem 4.1]{b4} & $[0,0,1,1],[0,0,2,2],\underline{2}[0,1,3,4].$\\

\hline $(\mathbb{Z}_{45},5,4)$-SDF & \cite[Lemma 2.8]{MaChang04} & $\underline{4}[0,2,5,12,23],\underline{4}[0,1,14,20,29]$,\\
&&$[0,4,4,-4,-4]$.\\

\hline $(\mathbb{Z}_{119},8,8)$-SDF & \cite[Lemma 3]{cfw2} & $ [20, 20,-20,-20, 29, 29,-29,-29]$,\\ & &
$\underline{4} [0, 1, 42, 28, 101, 97, 94, 114],$\\
& & $\underline{4} [0, 1, 12, 23, 41, 85, 104, 106]$, \\ &&
$\underline{4}[0, 2, 5, 17, 37, 47, 68, 76],$\\ & & $
\underline{4}  [0, 4, 10, 38, 54, 62, 86, 93].$\\\hline&&\\

\hline Type $6$ &&\\
\hline $(\mathbb{Z}_2,3,12)$-SDF & \cite[Theorem 3.1]{byw} & $[0,0,0],\underline{3}[1,1,0].$\\
\hline $(\mathbb{Z}_4,4,12)$-SDF & \cite[Case 2 in Section 2]{b4} & $\underline{3}[0,0,1,2],[0,0,0,1].$\\
\hline $(\mathbb{Z}_4,4,6p)$ & \cite[Construction A]{b4} & $\underline{(3p-9)/2}[0,1,2,3],\underline{(p-5)/2}[0,0,0,0],$\\
for an odd prime $p>3$ & & \underline{2}[0,1,1,1],\underline{3}[0,0,1,2],$[0,0,0,1],[0,0,0,2]$.\\
\hline $(\mathbb{Z}_8,4,6)$-SDF & \cite[Theorem 5.1]{b4} & $\underline{3}[0,1,3,5],[0,0,0,1].$\\\hline&&\\
\hline Type $7$ &&\\
\hline $(\mathbb{Z}_6,7,56)$-SDF & \cite[Theorem 3.4]{byw} &$\underline{7}[0,1,2,3,4,5,5],[0,0,0,0,0,0,0].$\\
\hline
\end{tabular}
\caption{SDFs of type $d$ in the literature}\label{tab:3}
\end{center}
\end{table}

There are some more papers where SDFs of specific type have been crucial. For example see \cite{b19,m}. Note, however, that the blocks of the SDFs in \cite{b19} do not have constant size.

On the other hand, even though we give the definition of strong difference families in general group,  only abelian SDFs are considered in this paper. We notice that non-abelian SDFs, especially the dihedral ones, have been investigated by Buratti and Gionfriddo in \cite{bg}. By using a $(D_{12}, 9, 6)$-SDF of type $6$, they constructed an infinite class of $9$-GDDs of type $12^u$. They also pointed out that via a $(D_6, 9, 12)$-SDF, Abel, Bluskov and Greig \cite{abg} constructed an infinite classes of $9$-GDDs of type $6^u$. Since in this paper we deal mainly with SDFs of type $d$ for $d\in\{2,4\}$, one possible future direction of work is to investigate carefully other values of $d$. In particular we hope to find constructions like the ones of Propositions \ref{pat2} and \ref{pat4} for $d\in\{3,6\}$ considering both abelian and non-abelian cases.

Finally, there should be no obstacle to extending Proposition \ref{prop:basic} to higher index $\lambda$. That is another further direction.

\section*{Aknowledgments}

The authors express their gratitude to the anonymous referees for their detailed and constructive comments which are very helpful to the improvement of the paper. Thank one of the referees for pointing out the dihedral SDFs used in reference \cite{bg}. The authors thank Professor Marco Buratti of Universit\`a di Perugia for his many valuable comments. Research of this paper was partially carried out while the second author was visiting Beijing Jiaotong University. He expresses his sincere thanks to the 111 Project of China (grant number B16002) for financial support and to the Department of Mathematics of Beijing Jiaotong University for their kind hospitality.

\end{document}